\documentclass{article}

\newcommand{\D}{\displaystyle}
\newcommand{\Z}{{\mathbb Z}}

\newtheorem{lemma}{\bf Lemma}{}{\rmfamily}
{}{\rmfamily}
{\itshape}{\rmfamily}
\newtheorem{definition}{\bf Definition}{}{\rmfamily}
{}{\rmfamily}
{}{\rmfamily}
\newtheorem{remark}{\bf Remark}{}{\rmfamily}
\newtheorem{theorem}{\bf Theorem}{}{\rmfamily}
\usepackage{amsmath}
\usepackage{amsfonts}
\usepackage{amssymb}

\begin{document}

\title{Equivalences of $\Z _t \times \Z_2^2$-cocyclic
Hadamard matrices}

\author{
V. Alvarez, F. Gudiel, M. B. Guemes, 
K. J. Horadam, A. Rao
%
}

\maketitle

\begin{abstract}
One of the most promising structural approaches to resolving the
Hadamard Conjecture uses the family of cocyclic matrices over $\Z
_t \times \Z_2^2$.  Two types of equivalence relations for
classifying cocyclic matrices over $\Z _t \times \Z_2^2$ have been
found. Any cocyclic matrix equivalent by either of these relations
to a Hadamard matrix will also be Hadamard.

One type, based on algebraic relations between cocycles over any
finite group, has been known for some time. Recently, and independently, a second
type, based on four geometric relations between diagrammatic
visualisations of cocyclic matrices over $\Z _t \times \Z_2^2$,
has been found. Here we translate the algebraic equivalences to
diagrammatic equivalences and show one of the diagrammatic equivalences
cannot be obtained this way. This additional equivalence is shown to be the geometric translation of matrix transposition. 

\end{abstract}

\section{Introduction}
A Hadamard matrix of order $m$ is a square matrix $[h(i, j)]$ with
entries $h(i, j) = \pm 1$, $1 \leq i, j \leq  m$, whose row
vectors are pairwise orthogonal. A Hadamard matrix must have order
$1$, $2$ or a multiple of $4$, but no other restrictions on the
order of a Hadamard matrix are known, and the century-old Hadamard
Conjecture proposes that a Hadamard matrix exists for every $m
\equiv  0 \pmod 4$.

About 20 years ago, the use of cocycles and cocyclic matrices was
introduced by Horadam and de Launey \cite{HdL95} as a structural
approach to resolving the Hadamard Conjecture. Its advantages led to the
cocyclic Hadamard conjecture: that a cocyclic Hadamard matrix exists for every $m
\equiv  0 \pmod 4$.  The study and use of cocyclic matrices has expanded substantially since then, to include generalised Hadamard matrices \cite{Hor07,Hor10} and pairwise combinatorial designs \cite{deLF11}.

If $G$ is a group
and $C$ is an abelian group, a ($2$-dimensional, normalized)
cocycle $\psi$ is a mapping $\psi : G \times  G \to C$ satisfying
$\psi (1, 1) = \psi (g, 1) = \psi (1, g) = 1, ~ g \in G$ and the
cocycle equation:
\begin{equation}
 \psi (g, h) ~\psi (gh, k) ~=~ \psi (g, hk)~ \psi(h, k), ~~g, h, k \in G.
\end{equation}
The set of cocycles from $G$ to $C$ forms an abelian group $Z^2(G,
C)$ under pointwise multiplication. The simplest cocycles are the coboundaries
$\partial f$, defined for any function $f : G \rightarrow C$  by
$\partial f(g, h) = f(g)^{-1} f(h)^{-1} f(gh)$.

A cocycle may be represented by its matrix of values
\begin{equation}\label{coc}
 M_\psi = [ \psi (g,h) ]_{g,h \in G}
\end{equation}
once an indexing of the elements of $G$ has been chosen.

We set $C = \{ \pm 1 \} \cong \Z_2$ when searching for cocyclic
Hadamard matrices. A cocycle $\psi$ for which the cocyclic matrix $M_\psi$ is Hadamard is termed {\it orthogonal}. It is computationally easy to check whether $M_\psi$ is
a Hadamard matrix, as we only need to check whether the dot
product of the first row with each other row is $0$. This
computational cutdown is one motivation for using cocyclic matrices.

Most of the known constructions of Hadamard matrix families are cocyclic \cite[Ch. 6]{Hor07}.
Computationally, the most prolific indexing groups $G$ for
producing cocyclic Hadamard matrices appear to be the abelian
groups $\Z _t \times \Z_2^2$ and the dihedral groups $D_{4t}$,
where we may assume $t$ is odd. The $D_{4t}$ family, related to the
Ito type Hadamard matrices, has been investigated
by many researchers including the authors (see \cite{Hor07}).
The  $\Z _t \times \Z_2^2$ family, related to the  Williamson type Hadamard matrices,
has also been investigated by the authors \cite{AGG11a,BH95}, and  while exhaustive search often finds fewer Hadamard matrices in each order than for $D_{4t}$, abelian-ness makes the family computationally more tractable.

In parallel with the search for examples of Hadamard matrices in new orders, whether cocyclic or not, has been the attempt to classify them into equivalence classes. Hadamard equivalence of a $\{\pm 1\}$ matrix involves only permutation of rows or columns, and multiplication of a row or column by $-1$. While the transpose of a Hadamard matrix is a Hadamard matrix, transposition is not a Hadamard equivalence. The total number of Hadamard equivalence classes in small orders grows so rapidly that Orrick \cite{Orr08} uses a coarser $Q$-equivalence relation on Hadamard matrices which allows extra ``switching" operations and leads to a dramatic reduction in the number of classes.

The total number of equivalence classes of cocyclic Hadamard matrices over all index groups $G$ is studied by  \'{O} Cath\'{a}in and R\"{o}der \cite{OCR11} and calculated up to $m = 36$. An allied but distinct approach has been to identify equivalences of cocycles that preserve orthogonality. For the  $\Z _t \times \Z_2^2$ family, two different types of equivalence of cocycles, both of which preserve orthogonality, have been discovered independently.

The first of these is defined (see \cite{Hor07}) for any $G$ and $C$ by all compositions of a ``shift" action and two ``automorphism" actions.  (If $C =  \{ \pm 1 \}$ one of the automorphism actions is
trivial.) The resulting equivalence classes, called {\em bundles},  are already known by other names in different contexts. For example, if $f$ is a cryptographic function and $\psi = \partial f$ is a coboundary, the bundle corresponds to the Extended Affine (EA) equivalence class of $f$. Shift action is also studied separately, for applications to the search for self-dual codes \cite{Rao05} and, via shift representations, to classification of pairwise combinatorial designs \cite{FlEg2014}.

The second of these equivalences, independently introduced in \cite{AGG11a}, is specific to cocycles $\psi$ in $Z^2 := Z^2(\Z _t \times \Z_2^2,  \{ \pm 1 \})$ and arises from detailed investigation of a generating set of cocycles for $Z^2$.   Corresponding to the decomposition of $\psi$ as a product of generators there is a Hadamard product decomposition of $M_\psi$  into generator matrices. Geometric actions on these generator matrices lead to a concise diagrammatic representation of cocycles and geometric equivalences which is very useful for effective computation.

This paper relates and reconciles the two types of equivalence.

The paper is organized as follows. Section \ref{sec:background} describes the two types of equivalence. The group acting on cocycles is determined for each type; the two groups are not isomorphic. Section \ref{sec:results} gives our main results,
Theorems \ref{thm:translate} and \ref{thm:main}, translating shift
action and the remaining automorphism action into diagram actions,  relating the two groups of actions, and showing that the diagram action termed ``complement" has no algebraic analogue. In Section \ref{sec:complement} this diagram action is shown to be the transposing operation on $M_\psi$. We summarise and suggest further work.

\section{Background} \label{sec:background}
From now on we assume $C = \{ \pm 1 \}$, $G \cong \Z _t \times \Z_2^2$ with $t>1$
odd,  and $\psi \in Z^2$.  Denote the group of units of the ring $\Z_t$ by $\Z_t^*$. Let $G$ have presentation $$G = \langle x,u,v:\;
x^t=u^2=v^2=1,xu=ux,xv=vx,uv=vu\rangle ,$$ and ordering
$$(x^i,1)<(x^i,u)<(x^i,v)<(x^i,uv), \, 0 \leq i < t, \; 
(x^i,uv) < (x^{i+1},1), \, 0 \leq i <t-1\,.$$
We describe an orthogonality-preserving algebraic action on $\psi$ in the first subsection and an orthogonality-preserving geometric action on $\psi$ in the second.

\subsection{Bundle action on cocycles}

For any $a \in G$, the {\it shift} $\psi \cdot a$ of $\psi$ is the cocycle $(\psi \cdot  a)(g, h) = \psi (ag, h) \psi(a, h)^{-1}.$ It is orthogonal if $\psi$ is orthogonal. For any automorphism $\theta \in {\rm Aut}(G)$, the cocycle $\psi \circ (\theta \times \theta)$ is orthogonal if $\psi$ is. When the two actions are combined, the result is an action by the semidirect product $H = G \rtimes {\rm Aut}(G)$ called {\em bundle action} under which the orbit of $\psi$ is its {\em bundle}
\begin{equation} \label{eq:bundle} {\cal B}(\psi) = \{(\psi \cdot a)
\circ (\theta \times \theta): \; a \in G, \; \theta \in
\mbox{Aut}(G)\}.
\end{equation}
The group $H$ acting on $Z^2$ is $H = G \rtimes {\rm Aut}(G)$, where the
semidirect product is defined for $a,b \in G, ~\theta_1,\theta_2 \in {\rm
Aut}(G)$ by $a \theta_1 \circ b \theta _2 = a\theta_1^{-1}(b)\theta_1\theta_2$. See \cite[Ch. 8]{Hor07} for details. \\

The Hadamard equivalence operations on $M_\psi$
corresponding to shift and automorphism action can be easily described. $M_{\psi \cdot a}$
is Hadamard equivalent to $M_\psi $ by first permuting the rows of
$M_\psi$ with respect to the row index permutation $g \mapsto g' =
ag, ~g \in G$, obtaining $M' = [\psi (ag, h)]_{g,h\in G}$. The
first row of $M'$ is the $a^{th}$ row of $M_\psi$. Then obtain
$M_{\psi \cdot a}$ from $M'$ by multiplying every row of $M'$
point-wise by its first row, or, equivalently, by multiplying
every column of $M'$ by its first entry. $M_{\psi \circ (\theta
\times \theta)}$ is Hadamard equivalent to $M_\psi$ by permuting
rows and columns under $\theta$.

We complete this subsection by identifying the group $H = G \rtimes
{\rm Aut}(G)$ which partitions cocycles into bundles
(\ref{eq:bundle}).

\begin{theorem} \label{thm:bundle group}
The group $H$ defined by bundle action on $Z^2$ is $H \cong [\Z_t \rtimes \Z_t^*] \times  [\Z_2^2 \rtimes S_3]$. Therefore the order of $H$ is $\,24\, t\, \phi(t)$,\,
where $\phi$ is the Euler function.

A generating set for $H$ is $\{x, u, v, h_r, r \in \Z _t^*, h_{23},  h_{243} \}$, where $x, u$ and $v$ are shift actions and $h_{23}: x \mapsto x, u \mapsto v, v \mapsto
u$; $h_{243}: x \mapsto x, u \mapsto uv, v \mapsto u$ and $h_{r}:
x \mapsto x^r, u \mapsto u, v \mapsto v$ are automorphism actions. 
\end{theorem}
{\it Proof.} Since $t$ is odd, ${\rm Aut}(\Z _t \times \Z_2^2) \cong {\rm Aut}(\Z_t)
\times {\rm Aut}(\Z_2^2) \cong \Z _t^* \times S_3$.
Under the identification $1 \leftrightarrow 1, u \leftrightarrow
2, v \leftrightarrow 3, uv \leftrightarrow 4$, ${\rm Aut}(\Z_2^2)$
is the subgroup of $S_4$ which fixes $1$. Then
$\{{\rm Id}\} \times {\rm Aut}(\Z_2^2)$ is generated by $h_{23}$
and $h_{243}$. Thus $H = [\Z _t \times \Z_2^2] \rtimes [\Z _t^*
\times S_3]$, with the listed generating set. Since $h_{23}(x) =
h_{243}(x) = x$, $\Z_t$ commutes with $S_3$ and since $h_r(u) =
u,\, h_r(v) = v$,  $\Z_2^2$ commutes with $\Z_t^*$. Hence $H \cong
[\Z_t \rtimes \Z_t^*] \times [\Z_2^2 \rtimes S_3]$. \hfill
$\square$

\begin{remark}\label{nota1} In terms of the Coxeter presentation of $S_n$, if $\sigma _i$ denotes the transposition $(i\; i+1)$, $S_n = \langle \sigma _i: \; \sigma _i^2=(\sigma _i \sigma _{i+1})^3=1, 1 \leq i \leq n-1 \rangle$ and $S_4=\langle \sigma _1, \sigma_2, \sigma_3: \; \sigma _i^2=(\sigma _i \sigma _{i+1})^3=1, 1 \leq i \leq 3 \rangle > \langle \sigma_2,\sigma_3 \rangle \cong S_3$, so that  in Theorem \ref{thm:bundle group}, $h_{23}=\sigma_2$ and $h_{243}=\sigma_3\sigma_2$. \end{remark}

\subsection{Geometric action on cocycle diagrams}

The group of cocycles $Z^2$ has a generating set ${\cal Z} = \{\partial _1, \ldots,
\partial _{4t}, \beta_1, \beta_2, \kappa\}$ consisting of $4t$
coboundaries $\partial _i := \partial \delta_i$, where $\delta_i$ is the Kronecker delta function of the $i^{th}$-element in $G$ in the given ordering, and three representative cocycles
$\beta_1, \beta_2, \kappa$, all of which are explicitly described in \cite{AAFR08,AGG11a}. Every 2-cocycle over $G$ admits a (non unique) representation as a product of the generators in
${\cal Z}$. The identity of $Z^2$ is the trivial cocycle  ${\bf 1}$ for which $M_{\bf 1} = J_{4t}$ is the all-ones matrix. All orthogonal cocycles known so far (cf. \cite{BH95,AGG11a}) contain the factor $\rho ={\beta_1}{\beta_2}{\kappa}$,  where
\begin{equation} \label{eq:M_rho}
M_\rho = J_t \otimes
 \left(\begin{array}{rrrr}1&1&1&1\\1&-1&1&-1 \\1&-1&-1&1\\1&1&-1&-1
\end{array}\right)\,,
\end{equation}
and $J_t$ denotes the $t \times t$ matrix all of 1s.   It is conjectured this must always be true \cite[Research Problem 37]{Hor07}. For the remainder of the paper, we assume that we work with cocycles of this type. That is, $\psi =
{\partial_1}^{\epsilon_1} \ldots
{\partial_{4t}}^{\epsilon_{4t}}\, \rho$, $\epsilon_i \in \{0, 1\}$.

In \cite{AGG11a},  a more concise notation to describe $\psi = {\partial_{d_1}} \ldots {\partial_{d_k}}\, \rho$ is introduced, which allows one to determine if $\psi$ is orthogonal much
more easily. Partition the set $\{d_1,
\ldots, d_k\}$ according to the equivalence classes modulo 4, in the
class order $2, 3, 0, 1$ and in descending order within each
class. We will denote this ordered set of coboundaries
\begin{equation}
\{ {\bf {c_2,c_3,c_4,c_1}} \} = \{ \{d_{2+4 j_2}\},\{d_{3+4 j_3}\},
\{d_{4 j_4}\}, \{d_{1+4 j_1}\} \}.
\end{equation}
For example, for $t = 7$, the cocycle $\; \psi = {\partial_{4}} {\partial_{6}}
{\partial_{9}}{\partial_{10}}{\partial_{11}}{\partial_{12}}
{\partial_{14}}{\partial_{20}}{\partial_{21}}{\partial_{25}}\,\rho\;$
is orthogonal, and is represented as
\begin{equation} \label{ex:t=7}
\{ {\bf {c_2,c_3,c_4,c_1}} \} = \{ \{14, 10, 6\}, \{11\}, \{20, 12, 4 \}, \{25, 21, 9 \} \}.
\end{equation}
Alternatively, we can write all the integers $1, \dots , 4t$, by
equivalence classes modulo 4, in descending order, as the rows of a $4 \times
t$ matrix (treated as a cylinder, i.e. left and right edges are
identified) and mark out only the entries occurring in $\{d_1, \ldots,
d_k\}$.
\begin{definition} \cite{AGG11a}
The {\em diagram} of $ \psi = {\partial_{d_1}} \ldots
{\partial_{d_k}} \, \rho$ is a $4 \times t$ matrix A, such that
$a_{ij} = \times$  if  $ 4t - 4(j-1)-3 + i \mod  4t \in \{ d_1,
\ldots , d_k  \}$ and $a_{ij} = -$ elsewhere.
\end{definition}

The diagram for the example in (\ref{ex:t=7}) above is
\begin{equation} \label{eq:diagram ex}
A = \left|
    \begin{array}{ccccccc}
        - & - & - & \times & \times & \times & - \\
       - & - & - & -  & \times & - & - \\
       - & - & \times & - & \times &  - & \times  \\
       \times & \times & - & - & \times & - & -\\
    \end{array}
  \right|
\end{equation}

We now list the four types of orthogonality-preserving operations on $\psi$ described in
\cite{AGG11a}. We adopt the notation $[m]_n$ for $m \mod n$ for brevity.

\begin{definition}
Let $\{ {\bf {c_2,c_3,c_4,c_1}} \}$ be a set of coboundaries. Denote the columns of its diagram $A$ by $ ({\cal{C}}_{t-1},\cdots ,{\cal{C}}_0)$.  Let $ {\bf {c_j + k }}$ denote the set of coboundaries obtained by adding $k$ to each element of ${\bf {c_j}}$ modulo $4t$.
\begin{enumerate}
\item The {\em complement} $ {\textsc{C}_2} (\{ {\bf {c_2},\bf{c_3,c_4,c_1}} \} )$ of this set is the set $\{ {\bf {\overline{c_2},c_3,c_4,c_1}} \} $
where ${\bf {\overline{c_2}}}$ is complement of ${\bf {c_2}}$ in the equivalence class 2 modulo 4.
\item Six elementary {\em swapping} operations are possible on this set: $s_{12}, s_{13}, s_{14}$ (see \cite{AGG11a}) and
\begin{itemize}
\item $s_{23} ( \{ {\bf {c_2,c_3,c_4,c_1}} \} ) = \{ {\bf {c_3 -1,c_2+1,c_4,c_1}} \}. $
\item $s_{24} ( \{ {\bf {c_2,c_3,c_4,c_1}} \} ) = \{ {\bf {c_4-2,c_3,c_2+2,c_1}} \}. $
\item $s_{34} ( \{ {\bf {c_2,c_3,c_4,c_1}} \} ) = \{ {\bf {c_2,c_4-1,c_3+1,c_1}} \}. $
\end{itemize}
\item The
{\em $i$-rotation} $\textsc{T}_i ( \{
{\bf {c_2,c_3,c_4,c_1}} \})$, $0 \leq i \leq t-1$,  of this set is the set
$$\{ {\bf c_2  -4i,  c_3 -4i, c_4 -4i, c_1 -4i} \}. $$
\item The {\em $r$-th dilatation} $\textsc{V}_r (\{ {\bf {c_2,c_3,c_4,c_1}} \})$, for $r \in \Z_t^*$, is the set with diagram $\textsc{V}_r(A)$, where
$\textsc{V}_r({\cal C}_{j}) = {\cal C}_{[jr]_t}, ~ 0 \leq j \leq t-1$.
\end{enumerate}
\end{definition}

Clearly the order of $\textsc{C}_2$  is 2 and $\langle \textsc{C}_2 \rangle \cong \Z_2$. The swappings each have order 2 and generate a group $\cong S_4$ which, in terms of a Coxeter presentation (Remark \ref{nota1}), is generated by  $\sigma_1 = s_{23}, ~\sigma_2 = s_{34}$ and $\sigma_3 = s_{14}$.  The rotations are generated by $\textsc{T}_1$ so $\langle \textsc{T}_1 \rangle \cong \Z_t$; and $\langle \textsc{V}_r, r \in  \Z_t^* \rangle \cong  \Z_t^*$.

In terms of diagrams, $\textsc{C}_2$  complements the first row of $A$; $s_{ij}$  swaps rows corresponding to ${\bf c_i}$ and ${\bf c_j}$; $\textsc{T}_i$ cyclically shifts
columns $i$ places to the right; and $\textsc{V}_r$ permutes columns according to multiplication of column index by the invertible element $r$ (so ${\cal C}_0$ is always fixed).

For instance, if $A$ is the diagram in (\ref{eq:diagram ex}),
{\small
$$
\begin{array}{ll}
 {\textsc{C}_2}(A) =
  \left|
    \begin{array}{ccccccc}
       \times & \times  & \times &  - & -  & -  & \times \\
       -& - & - &  - &  \times & - & - \\
       -& - &  \times & - & \times & - &  \times  \\
      \times & \times &  - &  - & \times  & - & -  \\
    \end{array}
  \right|,
&
\textsc{T}_2(A)  =  \left|
    \begin{array}{ccccccc}
         \times & - & - & - & - & \times & \times  \\
       -&  - & -  & - & -   & - & \times\\
       - &  \times & - &  - & \times  & - & \times \\
      - &  - & \times  & \times & -  & - & \times  \\
    \end{array}
  \right| , \\
  & \\
s_{23}(A) =   \left|
\begin{array}{ccccccc}
       -& - & - & - & \times  & - & -  \\
              -&  - & - &  \times & \times & \times  & - \\
       -& - &  \times & - & \times &  - & \times  \\
      \times & \times & - &  - & \times  & - & -  \\
    \end{array}  \right| ,
    &
 \textsc{V}_2(A) =   \left|
    \begin{array}{ccccccc}
           \times & - & \times  & - & \times  & -   & -\\
       -&  - & \times  & - & -   & -  & - \\
       - & -  & \times &  - & - & \times  &   \times \\
      - &  \times  & \times & \times  & -   & - & -  \\
    \end{array}
  \right| .
\end{array}
$$
}

It is possible to identify the action of $\textsc{C}_2$ on coboundaries directly.

\begin{lemma} \label{lem:C-action}
$ \textsc{C}_2( {\partial_{d_1}} \ldots {\partial_{d_k}}) =  {\partial_{d_1}} \ldots {\partial_{d_k}}~\prod _{i=0}^{t-1}\,{\partial_{2+4i}}.$
\end{lemma}
{\it Proof.}
If $\psi = {\bf 1}$  is the trivial coboundary in $Z^2$, with $M_{\bf 1} = J_{4t}$ the all-ones matrix, then it has  $\{ {\bf {c_2,c_3,c_4,c_1}} \} = \{ \emptyset, \emptyset, \emptyset, \emptyset\}$ so
$\textsc{C}_2({\bf 1})=\prod _{i=0}^{t-1}{\partial_{2+4i}}$. By simple inspection, it may be checked that
\begin{equation}\label{eq:C-action}
\textsc{C}_2(J_{4t})=\prod _{i=0}^{t-1}M_{\partial_{2+4i}}=J_t \otimes  \left(\begin{array}{rrrr}1&1&1&1\\1&1&-1&-1 \\1&-1&1&-1\\1&-1&-1&1
\end{array}\right).
\end{equation}
The result follows immediately.  \hfill $\square$

\

We complete this subsection by identifying the group $H'$ generated by the diagrammatic
operations above.

\begin{theorem} \label{thm:diagram group}
The group $H'$ defined by diagrammatic action on $Z^2$ is $H' \cong [\Z_t \rtimes \Z_t^*] \times S_4 \times \Z_2$. Therefore the order of $H'$ is $\,48 \, t\, \phi(t)$. A generating set for $H'$ is $\{ \textsc{T}_1,\, \textsc{V}_r, r \in  \Z_t^*,~ s_{14}, s_{23}, s_{34},\, \textsc{C}_2 \}$.
\end{theorem}
{\it Proof.} It is shown in \cite{AGG11a} that the
complement and swapping operations commute with each other and
with rotations and dilatations, but that  rotation and dilatation do not commute.
The composition $\textsc{V}_r^{-1}  \textsc{T}_1  \textsc{V}_r$ acts on column $[j]_t$ of $A$ to give column $[(jr-1)r^{-1}]_t=[j-r^{-1}]_t$, so $\textsc{V}_r^{-1}  \textsc{T}_1  \textsc{V}_r=\textsc{T}_{r^{-1}}$.  Define a homomorphism $\mu : \Z_t^* \rightarrow \mbox{Aut}(\Z_t)$ by $\mu (\textsc{V}_r)(\textsc{T}_1)=\textsc{T}_{r^{-1}}$. Consequently, $\langle \textsc{T}_1, \textsc{V}_r,   r \in  \Z_t^* \rangle \cong \Z_t \rtimes _\mu Z_t^*$.

Swapping permutes rows while rotations and
dilatations permute columns, so swapping is not in the subgroup of
$H'$ generated by rotations and dilatations. All combinations of swapping,
rotation and dilatation preserve the total number of coboundaries
but complementation does not, so complementation is not in the
subgroup of $H'$ generated by rotations, swappings and
dilatations. \hfill $\square$

\section{Bundle actions as Diagram actions} \label{sec:results}

In this section we express the bundle actions on $Z^2$ in terms of the diagrammatic
operations and identify the role of the diagrammatic action $\textsc{C}_2$. Subsection \ref{ssec:proof} is given to proving the following
theorem.

\begin{theorem} \label{thm:translate}
\begin{enumerate}
\item The shift actions by $x, u$ and $v$, respectively, on $\psi$, are the diagrammatic actions
$\textsc{T}_1,\,s_{12} s_{34}$ and $s_{13} s_{24}$, respectively.

\item The automorphism actions by $h_r, h_{23}$ and $h_{243}$,
respectively, on $\psi$, are the diagrammatic actions $\textsc{V}_{r^{-1}}, \, \textsc{C}_2
s_{23}$ and $s_{234} := s_{23}s_{24}$, respectively.
\end{enumerate}
\end{theorem}

From Theorem \ref{thm:translate} we obtain our main result.

\begin{theorem} \label{thm:main}
Bundle action by $H$ on $\Z_t \times \Z_2^2$-cocyclic matrices
corresponds to diagrammatic action by the subgroup $$H^* = \langle
\textsc{T}_1,\,\textsc{V}_{r^{-1}},\,s_{12} s_{34},\,s_{13} s_{24},\,\textsc{C}_2
s_{23},\,s_{23}s_{24}\rangle \cong (\Z_t \rtimes \Z_t^*) \times
S_4$$ of index 2 in $H'$. The operation $\textsc{C}_2$ is not in $H^*$.
\end{theorem}
{\it Proof.} Define a homomorphism $\alpha : H \rightarrowtail H'$
by $x \mapsto \textsc{T}_1, h_r \mapsto \textsc{V}_{r^{-1}},  u \mapsto s_{12}
s_{34}, v \mapsto s_{13} s_{24}, h_u \mapsto \textsc{C}_2 s_{23}$ and $h_v
\mapsto s_{23}s_{24}$.
By Theorem
\ref{thm:diagram group} and Theorem \ref{thm:translate}, $\alpha(\langle x, h_r, r \in \Z_t^*\rangle) = \langle \textsc{T}_1, \textsc{V}_{r^{-1}}, r \in \Z_t^*\rangle \cong  \Z_t \rtimes \Z_t^* $ is an isomorphism.

Let $CS_4$ be the subgroup of $H'$ isomorphic to $S_4$ which is
generated by the 6 order-2 elements $\textsc{C}_2 s_{ij}$ (i.e. compose
every transposition $s_{ij}$ with the complement $\textsc{C}_2$; they
commute so order doesn't matter). Products corresponding to even
permutations in $S_4$ will appear unchanged, while those
corresponding to odd permutations in $S_4$ will be multiplied by
$\textsc{C}_2$. Then, from Theorem \ref{thm:bundle group} and Theorem
\ref{thm:translate}, $\alpha(\Z_2^2 \rtimes S_3)$ is generated by
$\textsc{C}_2 s_{12} \textsc{C}_2 s_{34} = s_{12} s_{34}$ and $\textsc{C}_2 s_{13} \textsc{C}_2 s_{24}
=s_{13} s_{24}$ (shift action, isomorphic to $\Z_2^2$), and $\textsc{C}_2
s_{23}$ and $\textsc{C}_2 s_{23} \textsc{C}_2 s_{24} = s_{23}s_{24}$ (automorphism
action, isomorphic to $S_3$). Direct calculation shows that
$\alpha$ maps $\Z_2^2 \rtimes S_3$ onto $CS_4$, so $\alpha$ is an isomorphism.
Thus $H^* \cong (\Z_t \rtimes \Z_t^*)
\times S_4$, and $\alpha(H)$ does not contain $\textsc{C}_2$. \hfill
$\square$


\subsection{Proof of Theorem \ref{thm:translate}} \label{ssec:proof}

Every cocyclic matrix $M_\psi$ admits a
decomposition as the Hadamard (pointwise) product of the cocyclic
matrices corresponding to the generators. That is, $M_\psi =
M_{\partial_1}^{\epsilon_1} \ldots
M_{\partial_{4t}}^{\epsilon_{4t}}\, M_\rho$, $\epsilon_i \in \{0, 1\}$.

Each matrix $M_{\partial_i}$ is symmetric. Without loss of generality we negate the $i^{th}$ row and $i^{th}$
column of $M_{\partial_i}$. These Hadamard equivalent 
matrices, denoted $M_{i}$, have a very particular form (see
\cite{AAFR08} for details). Each $M_{i}$ is a $4 \times
4$-block back diagonal square matrix of order $4t$. The first
block row has a $4 \times 4$ matrix $A_{[i]_4}$ as the $\lceil
\frac{i}{4} \rceil ^{th}$ block and $4 \times 4$ all-1s blocks in
the other $t-1$ positions. The remaining block rows are obtained
by successively back-cycling the first.

The $4\times 4$-blocks $A_{[i]_4}$ depend on the equivalence class of $i$
modulo $4$, as follows. Let
$R = \left(
\begin{array}{rr} 0 & 1 \\ 1 & 0 \end{array}\right)$,  $D = \left(
\begin{array}{rr} -1 & 1\\ 1 &-1\end{array}\right)$ so  $DR =  \left(
\begin{array}{rr} 1 & -1\\ -1 &1\end{array}\right)$.
Then,  adopting the notation blank for 1 and $-$ for $-1$, for brevity,
$A_0=\left( \begin{array}{rr} &DR\\DR& \end{array}\right),
A_1=\left( \begin{array}{rr} D&\\ &D
 \end{array}\right),  A_2=\left(
\begin{array}{rr}DR& \\ &DR \end{array}\right),
A_3=\left( \begin{array}{rr}  &D\\D&\end{array}\right).$\\

It may be checked without difficulty that bundle action by each of $x, u, v, h_r$ and $h_{243}$ leaves $M_\rho$ invariant. Only action by $h_{23}$ alters $M_\rho$.  In terms of identifying diagram actions, it does not matter whether we work with $M_{\partial_i}$ or $M_i$ so we use the latter. We determine each bundle action on $M_i$ in the subsections below, concluding with the action of $h_{23}$ on $M_\rho$.

\subsubsection{Shift action of $x$}

First, we change the order of the elements in the group to $g'=xg$, obtaining
$$(x,1) < (x,u) < \dots < (x^{t-1},uv)<(1,1)< \dots < (1,uv) $$
that is, we put the  first block of $4$ elements at the end of the list.

For an individual coboundary $\partial_i$, the reordering takes
the first four rows to the last four, moving the other rows
upwards. Now the blocks $A_{[i]_4}$ start from the $\lceil
\frac{i}{4} \rceil -4 ^{th}$-column,  the negated row is the $i-4
^{th}$ row, and the negated column is still the $i^{th}$ column.
Next we perform the pointwise product of the first row and the
others. This first row (the former $5^{th}$) has two negative
entries, at positions $i$ and $i-4$, so we have to negate these
columns, getting the coboundary $\partial_{i-4}$.

So, the  action of $x$ on the cocyclic matrix is the
$1$-rotation $\textsc{T}_1$ on it.

\subsubsection{Shift action of $u$}

First, we change the order of the elements in the group to
$g'=ug$, obtaining
$$(x^i,u) < (x^i,1) < (x^i, uv) < (x^i, v),\,
0 \leq i < t, ~(x^i,v) < (x^{i+1},u),\, 0 \leq i <t-1,$$ that is,
we reorder every block of $4$ elements by means of the permutation
$\sigma = (12) (34)$.

For an individual coboundary $\partial_i$, the reordering permutes
rows in the same way. This permutation transforms the blocks
$A_{[i]_4}$ in the same way, under $(A_1 A_2)(A_3 A_0)$, the
negated row is the $\sigma(i)^{th}$ and the negated column is the
$i^{th}$. The  first row (the former $2^{nd}$) has two negative
entries, at positions $i$ and $\sigma(i)$. After negating these
columns, we  get the coboundary $\partial_{\sigma(i)}$.

So, the  action of $u$ on the cocyclic matrix is the
composition of swappings $s_{21} s_{34}$.

\subsubsection{Shift action of $v$}

First, we change the order of the elements in the group to $g'=vg$, obtaining
$$(x^i,v) < (x^i,uv) < (x^i, 1) < (x^i, u),\,
0 \leq i < t, ~(x^i,u) < (x^{i+1},v)\, 0 \leq i <t-1,$$ that is,
we reorder every block of $4$ elements by means of the permutation
$\sigma' = (13) (24)$.

For an individual coboundary $\partial_i$, the reordering permutes
rows in the same way. This permutation transforms the blocks
$A_{[i]_4}$ in the same way, under $(A_1 A_3)(A_2 A_0)$, the
negated row is the $\sigma'(i)^{th}$ and the negated column is the
$i^{th}$. The  first row (the former $3^{rd}$) has two negative
entries, at positions $i$ and $\sigma'(i)$. After negating these
columns, we  get the coboundary $\partial_{\sigma'(i)}$.

So, the  action of $v$ on the cocyclic matrix is the
composition of swappings $s_{13} s_{24}$.

\subsubsection{Automorphism action of $h_r$}

A straightforward algebraic calculation shows that $h_r(\partial
_k) = \textsc{V}_{r^{-1}}(\partial _k)$, for each
$k=x^{k_x}u^{k_u}v^{k_v}$. Set $\delta _{ij}=-1$ if $i=j$, and
$\delta _{ij}=1$ otherwise.

On one hand, $h_r(\partial
_k)(x^{i_x}u^{i_u}v^{i_v},~x^{j_x}u^{j_u}v^{j_v})$
\begin{eqnarray}
&=&
\partial _k( x^{r\cdot i_x \bmod t} \,u^{i_u} \,v^{i_v},~x^{r\cdot
j_x\bmod t}\,u^{j_u}\,v^{j_v}) \\\nonumber &=& \delta
_{x^{k_x}u^{k_u}v^{k_v},~x^{[r\cdot i_x]_ t}u^{i_u}v^{i_v}} ~
\delta _{x^{k_x}u^{k_u}v^{k_v},x^{[r\cdot
j_x]_t}u^{j_u}v^{j_v}}\\\nonumber & & \delta
_{x^{k_x}u^{k_u}v^{k_v},~x^{[r\cdot (i_x+j_x)]_
t}u^{[i_u+j_u]_2}v^{[i_v+j_v]_2}}.\label{delta1}
\end{eqnarray}

On the other hand, $\textsc{V}_{r^{-1}}(\partial
_k)(x^{i_x}u^{i_u}v^{i_v},~x^{j_x}u^{j_u}v^{j_v})$
\begin{eqnarray}
&=& \partial _{x^{[k_x\cdot r^{-1}]_t}\,u^{k_u}\,v^{k_v}}
(x^{i_x}\,u^{i_u}\,v^{i_v},~x^{j_x}\,u^{j_u}\,v^{j_v})
\\\nonumber &=&
\delta _{x^{[k_x\cdot
r^{-1}]_t}\,u^{k_u}\,v^{k_v},~x^{i_x}\,u^{i_u}\,v^{i_v}} ~ \delta
_{x^{[k_x\cdot
r^{-1}]_t}\,u^{k_u}\,v^{k_v},~x^{j_x}\,u^{j_u}\,v^{j_v}}
\\\nonumber & &
\delta _{x^{[k_x\cdot r^{-1}]_t}\,u^{k_u}\,v^{k_v},~x^{[i_x+j_x]_
t}\,u^{[i_u+j_u]_2}\,v^{[i_v+j_v]_2}}.\label{delta2}
\end{eqnarray}

A careful check, using the invertibility of $r$ in $\Z_t$, shows
these equations are equal term by term. Consequently, $h_r =
\textsc{V}_{r^{-1}}$, for all $r \in \Z_t^*$.

\subsubsection{Automorphism action of $h_{243}$}

The automorphism $h_{243}$ shifts cyclically to the right the
second, third and fourth positions of the elements in $G$, in each
block of $4$, leaving the first element unchanged. So the action
on the cocycles will be the same permutation of every second,
third and fourth rows and columns in every block of four.

For an individual coboundary $\partial_i$, this reordering
transforms the blocks $A_{[i]_4}$ in the same way, giving the
permutation $(A_2 A_3 A_0)$, and the negated row/column remains
unchanged if $[i]_4$ is $1$ and is interchanged cyclically between
cosets $2$, $3$ and $0$, so  we  get the coboundary $s_{234}
(\partial_i)$.

Hence, the action of  $h_{243}$ on any cocyclic matrix
gives us the operation $s_{234}$.

\subsubsection{Automorphism action of $h_{23}$}

The action of the automorphism $h_{23}$ on the cocyclic matrix
will be the permutation of second and third rows and columns in
every block of four.

For an individual coboundary $\partial_i$, this reordering
transforms the blocks $A_{[i]_4}$ in the same way, giving the
permutation  $(A_2 A_3)$, and the negated row/column remains
unchanged if $[i]_4$ is $0$ or $1$ and interchanged between cosets
$2$ and $3$, so  we  get the coboundary $s_{23} (\partial_i)$.

The action of this reordering on matrix $M_\rho$ applies the same
permutation to its rows and columns, so
the $4 \times 4$ blocks in (\ref{eq:M_rho}) become $$\left( \begin{array}{rrrr} 1&1&1&1 \\ 1&-1&-1&1\\
1&1&-1&-1 \\1&-1&1&-1 \end{array}\right).$$ This expression
coincides with the pointwise product of the $4 \times 4$ block in (\ref{eq:M_rho})  and the block $A_{2}$
with the second row and column negated, so the action of the
automorphism $h_{23}$ on $M_\rho$ gives us $M_\rho \cdot M_{\partial _2} \cdot
M_{\partial_ 6} \dots M_{\partial_{4t-2}}$, the product with all
coboundaries whose index is congruent to $2$ modulo $4$.
Hence,  by Lemma \ref{lem:C-action},
$h_{23}({\partial_{d_1}} \ldots {\partial_{d_k}}\, \rho) = s_{23}({\partial_{d_1}} \ldots {\partial_{d_k}})~(\prod _{i=0}^{t-1}\,{\partial_{2+4i}})\, \rho = \textsc{C}_2(s_{23}({\partial_{d_1}} \ldots {\partial_{d_k}}))\, \rho.$

Hence, the action of $h_{23}$ on any cocyclic matrix gives us the operation
$\textsc{C}_2 s_{23}$.

\section{Complement} \label{sec:complement}
Next we demonstrate that complementation corresponds to matrix transposition and gives the matrix of the transpose cocycle.
\begin{theorem} \label{thm:transposition}
The operation $\textsc{C}_2$ on $M_\psi$ coincides with transposition: $\textsc{C}_2(M_\psi)=(M_\psi) ^{\top} = M_{\psi^\top}$.
\end{theorem}
{\it Proof.} Consider $M_\psi = M_{\partial_{d_1}} \ldots M_{\partial_{d_k}}\, M_\rho$.
Since transposition commutes with pointwise products,
$M_\psi ^\top=M_{\partial_{d_1}} \ldots
M_{\partial_{d_k}}\, M_\rho^\top$.  By (\ref{eq:M_rho}) $$M_\rho^\top= J_t \otimes
 \left(\begin{array}{rrrr}1&1&1&1\\1&-1&-1&1 \\1&1&-1&-1\\1&-1&1&-1
\end{array}\right)= \left( J_t \otimes \left(\begin{array}{rrrr}1&1&1&1\\1&1&-1&-1 \\1&-1&1&-1\\1&-1&-1&1
\end{array}\right)\right)\, \cdot \, M_\rho \,.$$
By (\ref{eq:C-action}), $\displaystyle M_\psi ^\top = M_{\partial_{d_1}} \ldots
M_{\partial_{d_k}} \left(\prod _{i=0}^{t-1}M_{\partial_{2+4i}} \right) \, M_\rho = \textsc{C}_2(M_\psi)$, as claimed. Since $G = \Z _t \times \Z_2^2$ is abelian, the transpose $\psi^\top$ of $\psi$, with $\psi^\top(g, h) = \psi(h, g)$, is a cocycle \cite[(6.10)]{Hor07}, and $(M_\psi) ^{\top} = M_{\psi^\top}$. \hfill $\square$

\

In summary, we have shown that the diagrammatic operations which can be implemented for effective calculation of cocyclic Hadamard matrices over  $G = \Z _t \times \Z_2^2$, can all be interpreted as compositions of known algebraic equivalences, with the exception of complementation, which corresponds to matrix transposition. \'{O}Cath\'{a}in \cite{OCathainPhD} has used the algebraic equivalences together with transposition to determine classes of cocyclic matrices of order $4t$  over various $G$. He then checks any transposes lying in such a class to partition them  into Hadamard inequivalence classes. He coins the term {\it strong inequivalence} for Hadamard matrices $H$ and $H'$ for which $H'$ is not Hadamard equivalent to $H$ or to $H^\top$.  So, this approach using diagrammatic operations may be computationally effective.

It will also be interesting to investigate if diagrams and diagram operations can be found for cocycles over $G = D_{4t}$, and whether there are diagrammatic operations which correspond to Orrick's switching operations.

\end{document}